\pdfoutput=1
\documentclass[11pt]{article}

\usepackage[T1]{fontenc}
\usepackage[utf8]{inputenc}
\usepackage{lmodern}
\usepackage{amsmath,amssymb,amsthm,mathtools}
\usepackage{booktabs}
\usepackage{graphicx}
\usepackage{microtype}
\usepackage[round]{natbib}
\usepackage[colorlinks=true,allcolors=blue]{hyperref}
\usepackage{orcidlink}
\usepackage[margin=1in]{geometry}

\newtheorem{theorem}{Theorem}
\newtheorem{proposition}{Proposition}
\newtheorem{lemma}{Lemma}
\newtheorem{corollary}{Corollary}
\theoremstyle{remark}

\newcommand{\E}{\mathbb{E}}
\newcommand{\Pp}{\mathbb{P}}
\newcommand{\R}{\mathbb{R}}
\newcommand{\cF}{\mathcal{F}}
\newcommand{\cH}{\mathcal{H}}
\newcommand{\Fbar}{\overline{F}}
\newcommand{\gtwo}{\mathrm{g}_{2}}
\DeclareMathOperator{\Var}{Var}

\title{Moment-Orthogonal Quadratic-Envelope E-Processes\\
for Anytime-Valid Location Inference}
\author{Serhii Zabolotnii\,\orcidlink{0000-0003-0242-2234}\\[2pt]
\small Cherkasy State Business College, Cherkasy 18028, Ukraine\\
\small State Scientific Research Institute of Armament and Military Equipment\\
\small Testing and Certification, Cherkasy, Ukraine\\
\small Department of Cybernetics and Applied Mathematics, Uzhhorod University,\\
\small Uzhhorod 88000, Ukraine\\[2pt]
\small \texttt{zabolotnii.serhii@csbc.edu.ua}}
\date{31 October 2026}
\hypersetup{
  pdftitle={Moment-Orthogonal Quadratic-Envelope E-Processes for Anytime-Valid Location Inference},
  pdfauthor={Serhii Zabolotnii},
  pdfkeywords={confidence sequence, e-process, moment constraints, skewness}
}

\begin{document}
\maketitle

\begin{abstract}
E-process confidence sequences for unbounded means typically use a variance
bound and a first-order score, leaving skewness and kurtosis unused. We
construct a quadratic-envelope e-process from a conditionally centered vector
score with a predictable second-moment bound. Under moment information
specified independently of the tested stream, it stays positive on unbounded
support and gives exact anytime-valid coverage under predictable strategies
and fixed mixtures, without score clipping. For a location family with known
central moments through order four, the centered observation and its centered
square yield an exact confidence sequence whose local width shrinks, relative
to a first-order construction, by a factor determined by skewness and
kurtosis. At level 5\%, the largest null ever-crossing rate was 4.63\%; on
two zero-inflated-lognormal families median stopping times fell by
36.8--45.8\%. When the moments are calibrated on an independent reference
sample, a centering-robust factor keeps exact coverage at a closed-form width
cost: the median gain was 21.6--30.3\% for Gamma(2) and lognormal errors with
$10^6$ reference observations, and a loss of up to 9.9\% with $10^4$. On
public Fermi gamma-ray detector data with a declared stability margin, it
crossed at the current location in 1.8\% of 225 windows and stopped
30--35\% earlier than the calibrated first-order method.
\end{abstract}



\section{Introduction}

A confidence sequence covers its target simultaneously over time, so the
analyst may stop according to the data without invalidating the coverage
statement. Nonnegative supermartingales and Ville's inequality provide the
standard route to such guarantees \citep{ville1939,shafer2019,howard2021}.
Recent work has turned this route into a general language of safe,
anytime-valid inference \citep{ramdas2023}. Betting confidence sequences for
bounded observations use factors linear in $X-m$ \citep{waudby2024}; Catoni
supermartingales extend mean inference to unbounded observations when a
variance or fractional-moment bound is known \citep{catoni2012,wang2023}.

The restriction to first-order scores is structural.
For the class of finitely supported distributions on $\R$ with mean $m$, a
measurable, integrable function whose expectation vanishes under every member
of the class must be proportional to $x-m$.\footnote{Apply the condition to
two-point laws with mean $m$; measurability excludes pathological additive
solutions.} On the whole real line, requiring $1+\lambda(x-m)\geq0$ then leaves
only the trivial raw bet. Exact higher-order adaptation therefore requires a
null richer than the mean and scale information used by first-order and
Catoni-style constructions. Fan, Jiao and Wang constructed e-processes for
mean and variance constraints \citep{fan2025}. General theories now describe
e-variables for constraint-generated null classes and for conditionally
constrained sequential hypotheses \citep{larsson2025,clerico2026}. We study a
location null that additionally
specifies the third and fourth central moments. Under this null, the quadratic
envelope converts a moment-orthogonal score into an exact e-process on
unbounded support and gives an explicit local width constant.

Polynomial betting itself is not new. Legendre polynomial martingales have
recently been used to detect higher-order departures from exchangeability
\citep{szabadvary2026}. Nor is the scalar quadratic envelope new: it appears
in Catoni's influence-function construction. We combine that envelope with a
multivariate moment-orthogonal score and a predictable second-moment matrix
bound: the envelope supplies positivity, while the moment geometry supplies a
narrower confidence sequence on unbounded support. Here ``moment-enriched
null'' means that the candidate mean, variance, third central moment, and
fourth central moment are specified under the null; the last three are not
estimated from the tested stream. Our contribution is this combination and
its explicit location geometry. We do not claim that polynomial betting,
scalar quadratic envelopes, or e-processes under moment constraints are
themselves new. The growth analysis reduces local
optimization to the normal equations for the score covariance and alternative
mean, yielding the width constant $\sqrt{\gtwo}$ relative to the first-order
construction. The same geometry quantifies the fraction of local $\chi^2$
signal captured by the chosen score space, the sequential analogue of the
captured-fraction functional that governs estimation and testing with a
truncated likelihood expansion \citep{zabolotnii2026series}.

The general theorem needs a predictable matrix upper bound. The location
specialization additionally uses the exact variance to center the quadratic
score and the exact moment matrix to obtain $\sqrt{\gtwo}$. The main
simulations use these known null moments; a predictive plug-in variant, in
which they are estimated from the stream, lies outside the theorem.

The simulations with known moments corroborate the predicted width constant
and show earlier stopping on the tested skewed synthetic laws.

Oracle moments are then removed. If the moments are known only from an
independent reference sample, an exact confidence sequence still follows: a
finite-sample region for $(\sigma^2,\mu_3,\mu_4)$ under an eighth-moment
envelope, together with a factor that pays in advance for the remaining
centering error (Proposition~\ref{prop:robust}), keeps coverage at least
$1-\delta-\alpha_e$. Proposition~\ref{prop:calwidth} gives the resulting
width in closed form. The reference size $N$ limits the gain over a
first-order method calibrated on the same sample in two ways: uncertainty
about $\mu_3$ and $\mu_4$ caps it, and uncertainty about $\sigma^2$ erodes it
as the monitored length $n$ grows, with a loss of order $(n/N)^{1/2}$. The
gain vanishes once $n$ and $N$ are comparable. The reference need not have a
known location (Lemma~\ref{lem:region-unknown}), so a long record of the
process in a known-good state can serve as one. Section~\ref{sec:real-results}
applies the construction to photon energy channels of a gamma-ray detector.

The paper proceeds as follows. Section~\ref{sec:related} separates the
constraint frameworks relevant to the contribution.
Section~\ref{sec:factor} gives the e-factor and
its finite-sample validity. Section~\ref{sec:growth} derives the local growth
geometry and its remainder. Section~\ref{sec:location} specializes the result
to skewed location families, first with oracle moments and then, in
Section~\ref{sec:calibrated}, with moments calibrated on an independent
reference sample. Sections~\ref{sec:experiments} and
\ref{sec:results} describe the numerical study and its results.
Section~\ref{sec:discussion} states the scope that the method does and does not
cover.

\section{Relation to existing constraint frameworks}\label{sec:related}

Three kinds of constraints must be distinguished. \citet{larsson2025}
characterize e-variables for null classes generated by measurable
constraints. \citet{clerico2026} gives a sequential complete-class
result for hypotheses described by finitely many conditional constraints:
every e-process in that setting is pointwise dominated by a predictable
product of affine one-step e-variables. These results
concern restrictions on the null law. By contrast, \citet{saha2026} impose
structural constraints on the e-variable itself and study an
optimize-then-constrain principle. We use the first type of
information, namely conditional moment restrictions, and do not derive an
optimal e-variable over either general class.

Identification functions provide another route to sequential tests.
\citet{casgrain2022} construct supermartingales for elicitable and identifiable
functionals, allow predictable mixing, and study online convex optimization
strategies. Our narrower result gives a
positive polynomial factor on an unbounded support and works out the
moment-orthogonal location geometry and its width constant. The adaptive
strategies below are an implementation within the predictable-strategy
principle, not a claim that predictable mixing is new.

The closest specialized constructions delimit the remaining scope.
\citet{catoni2012} supplies the scalar positivity envelope, and
\citet{wang2023} adapt its influence function to confidence sequences.
\citet{fan2025} give exact joint mean--variance e-processes, while
\citet{szabadvary2026} uses polynomial martingales for exchangeability
testing. Positive-semidefinite matrix
supermartingales provide related matrix-valued concentration machinery
\citep{wangramdas2025}, but we do not identify the scalar factor here as a
special case of that framework.

The moment geometry itself comes from estimation theory. Kunchenko's
polynomial maximization method (PMM) estimates parameters from a finite set
of moments with stochastic polynomials, that is, weighted combinations of
centered powers of the observations \citep{kunchenko2002polynomial}. It has
been applied to signal parameters in non-Gaussian noise
\citep{kunchenko2003antenna,palahin2016joint}. For a location parameter, the
degree-two stochastic polynomial is the two-component score of
Section~\ref{sec:location}, and its asymptotic variance relative to the
sample mean is the coefficient $\gtwo$ of \eqref{eq:g2}. The same geometry
prices skewness in nonlinear system identification. A Volterra kernel
estimator built on the Gaussian-matched basis incurs an excess risk governed
by the third central moment of the input, which vanishes exactly when that
moment does; orthogonalizing the basis in the actual input law removes it
\citep{zabolotnii2026skew}. Likewise, $\gtwo=1$ exactly when $\gamma_3=0$.
We carry this
geometry from point estimation to anytime-valid inference: the e-process
uses the same score and moment matrix, but its guarantee is finite-sample and
does not rest on the asymptotics of the estimator.

\section{Quadratic-envelope e-processes}\label{sec:factor}

Let $(\cF_t)_{t\geq0}$ be a filtration. Predictable below means
$\cF_{t-1}$-measurable, so past-dependent estimates or bounds are allowed but
the current observation is not. At time $t$, let
$\phi_t\in\R^d$ be $\cF_t$-measurable and square-integrable and let
$\Fbar_t$ be a predictable positive semidefinite matrix with finite entries
almost surely. The
null class $\cH_0$ satisfies, for every
$P\in\cH_0$,
\begin{equation}\label{eq:null}
  \E_P(\phi_t\mid\cF_{t-1})=0,
  \qquad
  \E_P(\phi_t\phi_t^\top\mid\cF_{t-1})\preceq\Fbar_t.
\end{equation}
For any predictable vector $\lambda_t$ with finite entries almost surely, write
$h_t=\lambda_t^\top\phi_t$ and
$v_t=\lambda_t^\top\Fbar_t\lambda_t$, and assume $v_t<\infty$ almost
surely. Then $h_t^2$ is conditionally
integrable under \eqref{eq:null}. Define
\begin{equation}\label{eq:factor}
  E_t(\lambda_t)=\frac{q(h_t)}{1+v_t/2},
  \qquad q(u)=1+u+u^2/2.
\end{equation}

\begin{theorem}[Quadratic-envelope moment e-process]\label{thm:eprocess}
Under \eqref{eq:null}, if $\lambda_t$ is predictable with finite entries and
$v_t<\infty$ almost surely, then, for every $P\in\cH_0$,
$\E_P\{E_t(\lambda_t)\mid\cF_{t-1}\}\leq1$. Hence
$M_n=\prod_{t=1}^n E_t(\lambda_t)$ is a nonnegative supermartingale under
every $P\in\cH_0$, and
\[
  \sup_{P\in\cH_0}
  \Pp_P\!\left(\sup_{n\geq1}M_n\geq1/\alpha\right)\leq\alpha.
\]
If $M_n^{(k)}$ uses its own predictable strategy, then every finite or
countable mixture $\widetilde M_n=\sum_k w_kM_n^{(k)}$ with nonnegative,
data-independent weights summing to one satisfies the same bound.
\end{theorem}

\begin{proof}
Completing the square gives
$q(u)=\{(u+1)^2+1\}/2\geq1/2$, so the factor is strictly positive for every
real score. By predictability and \eqref{eq:null},
\[
 \E_P\{q(h_t)\mid\cF_{t-1}\}
 =1+\E_P(h_t\mid\cF_{t-1})+
   \tfrac12\E_P(h_t^2\mid\cF_{t-1})
 \leq1+v_t/2.
\]
The product is therefore a nonnegative supermartingale. Ville's inequality
gives the displayed probability bound. Linearity proves the finite-mixture
claim, and conditional monotone convergence proves it for a countable sum of
nonnegative components.
\end{proof}

The strategy need not be fixed in advance. For example, a regularized FTRL
update may choose $\lambda_t$ from
$\phi_1,\ldots,\phi_{t-1}$ and then use that value in the factor at time
$t$. Theorem~\ref{thm:eprocess} still applies because the update is
predictable. A cap on $v_t$ or regularization of the update can improve
numerical stability; neither is needed to make the quadratic envelope
positive. Using $\phi_t$ to choose $\lambda_t$ before evaluating the same
factor would violate predictability and is not covered.

The denominator in \eqref{eq:factor} compensates the quadratic envelope
exactly. Catoni's construction instead applies the bound
$1+v/2\leq\exp(v/2)$. This distinction slightly changes log-growth but not
the central conceptual requirement: validity depends on the moment constraint
in \eqref{eq:null}.

For each candidate $\theta$, let $\cH_0(\theta)$ denote the laws satisfying
\eqref{eq:null} for $\phi_t(\theta)$ and $\Fbar_t(\theta)$. Suppose that
$\theta\mapsto\phi_t(\theta)$ and $\Fbar_t(\theta)$ are jointly measurable,
and that the strategies are predictable for every fixed $\theta$. Assume also
that the displayed
supremum below is measurable. If $P\in\cH_0(\theta(P))$ at the true value,
inversion gives
\[
  C_n=\left\{\theta:
  \sup_{s\leq n}M_s(\theta)<1/\alpha\right\}.
\]
Theorem~\ref{thm:eprocess} implies
$\Pp_P(\theta(P)\in C_n\ \text{for all }n)\geq1-\alpha$.
No asymptotic argument enters this statement. Appendix~\ref{app:measurability}
separates pointwise coverage from the additional regularity needed to report
measurable scalar endpoints.

The oracle moment requirement can be transferred from an independent
calibration sample without treating a plug-in estimate as exact. Write
$\eta=(\sigma^2,\mu_3,\mu_4)$ and let
$\mathcal U=\mathcal U(D_{\rm cal})$ be a random confidence region based on
data independent of the tested stream.

\begin{corollary}[Independent calibration transfer]\label{cor:calibration}
Suppose $\Pp(\eta_0\in\mathcal U)\geq1-\delta$ and, on the event
$\{\eta_0\in\mathcal U\}$, conditional on $D_{\rm cal}$,
$M_n(m_0,\eta_0)$ is an e-process for the independent test stream. For
$0<\alpha_e<1$, define
\[
 C_n^{\rm cal}=
 \left\{m:\text{there exists }\eta\in\mathcal U\text{ such that }
   \sup_{s\leq n}M_s(m,\eta)<1/\alpha_e\right\}.
\]
Then the event $\{m_0\in C_n^{\rm cal}\text{ for every }n\}$ contains an
event of probability at least $1-\delta-\alpha_e$; in particular
\[
 \Pp\{m_0\in C_n^{\rm cal}\text{ for every }n\}
 \geq 1-\delta-\alpha_e
\]
whenever the left-hand event is measurable. With $\alpha_e=\alpha-\delta$ the
coverage is at least $1-\alpha$ whenever $0<\delta<\alpha$.
\end{corollary}

\begin{proof}
Let $G=\{\eta_0\in\mathcal U\}$ and
$H=\{\sup_{s\geq1}M_s(m_0,\eta_0)<1/\alpha_e\}$, both measurable. On $G\cap H$
the witness $\eta=\eta_0$ puts $m_0$ in $C_n^{\rm cal}$ for every $n$, so only
validity at the true pair $(m_0,\eta_0)$ is used and no simultaneous e-process
claim over the profiled values is needed. Conditional Ville inequality gives
$\Pp(G\cap H^c)\leq\alpha_e$, and a union bound with
$\Pp(G^c)\leq\delta$ gives $\Pp(G\cap H)\geq1-\delta-\alpha_e$.
\end{proof}

Corollary~\ref{cor:calibration} is a transfer principle, not a
distribution-free recipe for constructing $\mathcal U$. On an unbounded
support, finiteness of the fourth moment alone does not yield a nontrivial
finite-sample distribution-free upper confidence bound for that moment.
Useful regions therefore require additional tail or higher-moment
information, bounded calibration variables, a parametric restriction, or a
deterministic domain constraint. The obstruction is elementary: a finite
calibration sample can miss an arbitrarily rare component that carries an
arbitrarily large fourth moment.

Profiling over a continuous region $\mathcal U$ also has a computational
cost: the infimum over $\eta$ must be evaluated for every candidate. When the
calibration error enters only through the centering of the score, it can
instead be charged to the factor itself.

\begin{proposition}[Centering-robust factor]\label{prop:robust}
Let $q$, $h_t=\lambda_t^\top\phi_t$ and $v_t$ be as in
Theorem~\ref{thm:eprocess}, with $\phi_t$ adapted. Suppose that, instead of
\eqref{eq:null}, every $P\in\cH_0$ satisfies
\[
 \E_P(\phi_t\mid\cF_{t-1})\in\mathcal B_t,
 \qquad
 \E_P(\phi_t\phi_t^\top\mid\cF_{t-1})\preceq\Fbar_t,
\]
for a predictable set $\mathcal B_t\subset\R^d$. Let $\lambda_t$ be
predictable with finite entries and $v_t=\lambda_t^\top\Fbar_t\lambda_t<\infty$
almost surely, and let $\rho_t$ be predictable with
$\sup_{b\in\mathcal B_t}\lambda_t^\top b\leq\rho_t<1/2$. Then
\[
 E_t^{\rm rob}(\lambda_t)=\frac{q(h_t)-\rho_t}{1+v_t/2}
\]
is strictly positive and
$\E_P\{E_t^{\rm rob}(\lambda_t)\mid\cF_{t-1}\}\leq1$ for every
$P\in\cH_0$. Products and data-independent mixtures of such factors are
therefore e-processes, as in Theorem~\ref{thm:eprocess}.
\end{proposition}

\begin{proof}
Positivity follows from $q(h_t)\geq1/2>\rho_t$ and
$0\leq v_t<\infty$. Conditionally on
$\cF_{t-1}$, with $b_t=\E_P(\phi_t\mid\cF_{t-1})$,
\[
 \E_P\{q(h_t)\mid\cF_{t-1}\}-\rho_t
 =1+\lambda_t^\top b_t-\rho_t+\tfrac12\lambda_t^\top
   \E_P(\phi_t\phi_t^\top\mid\cF_{t-1})\lambda_t
 \leq1+v_t/2.
\]
\end{proof}

Proposition~\ref{prop:robust} is the exact counterpart of
\eqref{eq:misspec} below: a centering error $b_t$ that is known only to lie
in $\mathcal B_t$ is paid for in advance through $\rho_t$, and the second
moment matrix bound must include $b_tb_t^\top$. With $\mathcal B_t=\{0\}$
and $\rho_t=0$ it reduces to Theorem~\ref{thm:eprocess}.

\section{Growth geometry}\label{sec:growth}

We now fix one time point and take
$F=\E_{P_0}(\phi\phi^\top)\succ0$ to be the exact Gram matrix under a simple
reference null $P_0$. This exact $F$ is used for the growth and projection
identities; replacing it by a conservative upper bound preserves e-validity
but changes their numerical values. Let $Q$ be an alternative with
$\beta=\E_Q\phi$. The expected log-growth of
\eqref{eq:factor} is
\[
  G_Q(\lambda)=\E_Q\log q(\lambda^\top\phi)
  -\log\{1+\lambda^\top F\lambda/2\}.
\]

\begin{lemma}[Quantified quadratic approximation]\label{lem:remainder}
Let $h=\lambda^\top\phi$ and $v=\lambda^\top F\lambda$. If
$\E_Q|h|^3<\infty$, then
\begin{equation}\label{eq:remainder}
  G_Q(\lambda)=\lambda^\top\beta-v/2+R(\lambda),
  \qquad
  |R(\lambda)|\leq\frac{\E_Q|h|^3}{3}+\frac{v^2}{8}.
\end{equation}
\end{lemma}

\begin{proof}
For $r(u)=\log q(u)-u$,
$r'(u)=-u^2/\{2q(u)\}$. Since $q(u)\geq1/2$,
$|r'(u)|\leq u^2$, and integration from zero to $u$ yields
$|r(u)|\leq|u|^3/3$. Also
$0\leq x-\log(1+x)\leq x^2/2$ for $x\geq0$. Apply the latter inequality
with $x=v/2$ and combine the two bounds.
\end{proof}

The coefficient $1/2$ in the envelope has a convenient remainder bound, but
it is not singled out by second-order local growth.

\begin{proposition}[Local equivalence of quadratic envelopes]
\label{prop:quadratic-family}
For $c\geq1/4$, let
\[
 E_{t,c}=\frac{1+h_t+c h_t^2}{1+c v_t}.
\]
Under \eqref{eq:null}, $E_{t,c}$ is a nonnegative e-factor; it is strictly
positive when $c>1/4$. For the log-growth statement, assume $c>1/4$. Let
$\lambda=\epsilon a$ and consider a local path
$Q_\epsilon$ such that
\[
 \E_{Q_\epsilon}\phi=\epsilon b+o(\epsilon),
 \qquad
 \E_{Q_\epsilon}(\phi\phi^\top)=F+o(1).
\]
If, for some $\epsilon_0>0$,
$\sup_{|\epsilon|\leq\epsilon_0}\E_{Q_\epsilon}|a^\top\phi|^3<\infty$, then
\[
 \E_{Q_\epsilon}\log E_{t,c}
 =\epsilon^2\{a^\top b-a^\top Fa/2\}+o(\epsilon^2),
\]
independently of $c$. Thus $a^*=F^{-1}b$ and the leading optimal growth
$b^\top F^{-1}b/2$ are common to this quadratic family.
\end{proposition}

\begin{proof}
The discriminant of $1+h+c h^2$ is $1-4c\leq0$, which proves
nonnegativity; it is negative when $c>1/4$. Conditional centering and the
second-moment bound give
$\E(1+h_t+c h_t^2\mid\cF_{t-1})\leq1+c v_t$.
Finally,
\[
 \log(1+h+c h^2)=h+(c-1/2)h^2+O(|h|^3),
 \qquad
 \log(1+c v)=c v+O(v^2).
\]
Along the stated local path,
$\E h=\epsilon^2a^\top b+o(\epsilon^2)$,
$\E h^2=\epsilon^2a^\top Fa+o(\epsilon^2)$, and
$v=\epsilon^2a^\top Fa$. The terms involving $c$ cancel.
\end{proof}

Proposition~\ref{prop:quadratic-family} does not establish global optimality
among polynomial calibrators or among constrained e-variables. We retain
$c=1/2$ because it is strictly positive on the real line and gives the simple
cubic bound in Lemma~\ref{lem:remainder}.

If $F\succ0$, the quadratic term in \eqref{eq:remainder} has the unique maximizer
\begin{equation}\label{eq:normal}
  \lambda^*=F^{-1}\beta,
  \qquad
  G_{\mathrm{quad}}(\lambda^*)=\frac{1}{2}\beta^\top F^{-1}\beta=J/2.
\end{equation}
Thus the normal equations used for moment-orthogonal estimation also determine
the local growth-optimal score direction. The remainder states the region in
which that interpretation is accurate; it is not needed for e-validity.

The geometry is especially transparent when the score components are centered
and square-integrable under $P_0$. Let $z=dQ/dP_0-1$ in $L^2(P_0)$, and let
$\Pi z$ denote its orthogonal projection onto
$\operatorname{span}\{\phi_1,\ldots,\phi_d\}\subset L^2(P_0)$. Then
\begin{equation}\label{eq:kappa}
 J=\|\Pi z\|_2^2,
 \qquad
 \frac{\chi^2(Q\|P_0)-J}{2}
 =\left(1-\kappa\right)\frac{\chi^2(Q\|P_0)}{2},
 \quad
 \kappa=\frac{\|\Pi z\|_2^2}{\|z\|_2^2}.
\end{equation}
This identity is exact in the $\chi^2$ divergence. A corresponding Kullback--Leibler
statement requires a quadratic-mean differentiable local path; without that
condition, replacing $\chi^2/2$ by KL is not justified.

The same captured fraction governs fixed-sample procedures built from a
truncated likelihood expansion: there the captured fraction of the score is
the asymptotic relative efficiency and the local power of both estimation and
detection \citep{zabolotnii2026series}. Equation~\eqref{eq:kappa} is its
counterpart for sequential betting, where $\kappa$ is the share of the local
$\chi^2$ growth that the chosen score space retains.

\section{Skewness-adaptive location inference}\label{sec:location}

Let $X=m+\varepsilon$, where the centered error has known variance
$\sigma^2>0$, known third central moment $\mu_3$, and known finite fourth
central moment $\mu_4<\infty$. These three quantities are parts of the null
specification, not estimates from the tested stream. Exact $\sigma^2$ is
required to center the quadratic score; the more general matrix-bound
formulation of Theorem~\ref{thm:eprocess} does not by itself supply that
centering. Use
\[
 \phi_m(X)=
 \begin{pmatrix}X-m\\(X-m)^2-\sigma^2\end{pmatrix},
 \qquad
 F=
 \begin{pmatrix}
 \sigma^2&\mu_3\\
 \mu_3&\mu_4-\sigma^4
 \end{pmatrix}.
\]
Theorem~\ref{thm:eprocess} applies on an unbounded support. Under the alternative
$X=m+\delta+\varepsilon$, with the same centered error law as under the null,
$\beta(\delta)=(\delta,\delta^2)^\top$. Writing
$\gamma_3=\mu_3/\sigma^3$ and
$\gamma_4=\mu_4/\sigma^4-3$, define
\begin{equation}\label{eq:g2}
 \gtwo=1-\frac{\gamma_3^2}{2+\gamma_4}.
\end{equation}
This is the variance reduction coefficient of the degree-two polynomial
maximization method \citep{kunchenko2002polynomial}. Provided
$\Var(\varepsilon^2)>0$, positive definiteness of $F$ is equivalent to
$\gtwo>0$.

\begin{proposition}[Local width constant]\label{prop:width}
For the two-component score above,
\[
 \beta(\delta)^\top F^{-1}\beta(\delta)
 =\frac{\delta^2}{\sigma^2\gtwo}+O(\delta^3).
\]
For the comparison below, call two mixtures matched when they use the same
mixture weights, crossing level, and local boundary approximation after
score-specific normalization.

For the first-order score $X-m$, the corresponding leading term is
$\delta^2/\sigma^2$. Consider matched time-uniform mixtures whose local
half-widths $w_j(n)\to0$ satisfy
$n a_jw_j(n)^2=2b_n\{1+o(1)\}$ with the same boundary sequence $b_n$, where
$a_{\rm linear}=1/\sigma^2$ and
$a_{\rm moment}=1/(\sigma^2\gtwo)$. Then
\[
 \frac{\operatorname{width}_{\mathrm{moment}}(n)}
      {\operatorname{width}_{\mathrm{linear}}(n)}
 \longrightarrow\sqrt{\gtwo}.
\]
\end{proposition}

\begin{proof}
The $(1,1)$ entry of $F^{-1}$ is
$(\mu_4-\sigma^4)/\{\sigma^2(\mu_4-\sigma^4)-\mu_3^2\}
=1/(\sigma^2\gtwo)$. Since the second component of $\beta(\delta)$ is
$O(\delta^2)$, the first display follows. The stated matched-boundary
condition gives
$w_{\rm moment}/w_{\rm linear}\to
\sqrt{a_{\rm linear}/a_{\rm moment}}=\sqrt{\gtwo}$.
\end{proof}

The second score acts as a control variate. It has no first-order sensitivity
to a location shift, but its null correlation with $X-m$ removes a fraction
$1-\gtwo$ of the first-order score variance. For a symmetric law,
$\gamma_3=0$ and the asymptotic gain vanishes.

\subsection{Exact diagnostic under moment misspecification}
\label{sec:misspec}

The failure mode caused by incorrect moments can be written without an
asymptotic approximation. Suppose the implemented score and matrix are
$\widetilde\phi_t$ and $\widetilde F_t$, and define
\[
 b_t=\E(\widetilde\phi_t\mid\cF_{t-1}),\qquad
 S_t=\E(\widetilde\phi_t\widetilde\phi_t^\top\mid\cF_{t-1}).
\]
For a predictable $\lambda_t$, the implemented factor satisfies
\begin{equation}\label{eq:misspec}
 \E(\widetilde E_t\mid\cF_{t-1})
 =
 \frac{1+\lambda_t^\top b_t+
   \lambda_t^\top S_t\lambda_t/2}
 {1+\lambda_t^\top\widetilde F_t\lambda_t/2}.
\end{equation}
Thus centering error enters linearly through
$\lambda_t^\top b_t$. If $b_t=0$ and
$S_t\preceq\widetilde F_t$, conservative matrix inflation preserves
validity and loses growth. A nonzero $b_t$ or an underestimated matrix can
make the conditional mean exceed one. Because the error is repeated, a small
one-step excess can accumulate over a long horizon.

\subsection{Moments calibrated on an independent reference sample}
\label{sec:calibrated}

We now replace the oracle moments by a reference sample. The reference must
be external to the monitored data: estimated from the same observations, the
variance that centers the second score component removes the gain of the
degree-two polynomial location estimate altogether
\citep{zabolotnii2026rankone}. Suppose an
instrument is first run on a reference standard of known value, which yields
errors $\varepsilon_1',\ldots,\varepsilon_N'$ that are independent and
identically distributed, independent of the tested stream, and distributed as
the stream error $\varepsilon$, with known zero mean. The only tail
information assumed is an envelope $K_8\geq\E\varepsilon^8$. Split the
reference sample into a first part of size $N_A$ and a second part of size
$N_B=N-N_A$.

Let $Y_1,\ldots,Y_N$ be independent copies of a variable $Y$ with
$\Var(Y)\leq v$, fix a level $\delta$, and put
$\ell=\{2\log(1/\delta)/(Nv)\}^{1/2}$. Let $\psi$ be the influence function
of Section~\ref{sec:experiments}. Since
$e^{\psi(u)}\leq q(u)$ for every real $u$, the product
\[
 \prod_{i\leq N}\frac{e^{\psi\{\ell(Y_i-\E Y)\}}}{1+\ell^2v/2}
\]
has expectation at most one; the same holds with $-\psi$. By Markov's
inequality,
\[
 L=\inf\Bigl\{\theta:\sum_{i\le N}\psi\{\ell(Y_i-\theta)\}
   <N\log(1+\ell^2v/2)+\log(1/\delta)\Bigr\}
\]
is a lower confidence bound for $\E Y$ with error at most $\delta$, and the
mirror construction gives an upper bound. Because $\psi$ is nondecreasing,
both bounds are roots of monotone equations. We call them Catoni bounds
\citep{catoni2012}.

\begin{lemma}[Calibration region]\label{lem:region}
Let $K_8>0$ satisfy $\E\varepsilon^8\leq K_8$. Use the first split to compute an upper Catoni bound $\widehat\mu_4^{A}$ for
$\E\varepsilon^4$ with $v=K_8$. On the second split compute two-sided Catoni
intervals $[c-r,c+r]$ for $\sigma^2$ with $v=\widehat\mu_4^{A}$ and
$[c_3-r_3,c_3+r_3]$ for $\mu_3$ with $v=(\widehat\mu_4^{A}K_8)^{1/2}$, and an
upper bound $U_4$ for $\mu_4$ with $v=K_8$. If the four error levels sum to
$\delta$, then with probability at least $1-\delta$ all four statements hold
simultaneously.
\end{lemma}

\begin{proof}
$\Var(\varepsilon^4)\leq\E\varepsilon^8\leq K_8$ validates the first bound.
On its event, $\Var(\varepsilon^2)\leq\E\varepsilon^4\leq\widehat\mu_4^{A}$
and, by the Cauchy--Schwarz inequality,
$\Var(\varepsilon^3)\leq\E\varepsilon^6\leq(\E\varepsilon^4
\E\varepsilon^8)^{1/2}\leq(\widehat\mu_4^{A}K_8)^{1/2}$. These proxies are
fixed given the first split, which is independent of the second, so the
remaining bounds hold conditionally at their stated levels. A union bound
finishes the proof.
\end{proof}

A reference standard of known value is not always available, whereas a long
record of the process in a known-good state often is. The location of such a
record is unknown, but central moments do not depend on it.

\begin{lemma}[Reference of unknown location]\label{lem:region-unknown}
Let the reference record be $Y_i=\mu_{\mathrm{ref}}+\varepsilon_i'$ with
$\mu_{\mathrm{ref}}$ unknown and the $\varepsilon_i'$ as above. Split it into
parts of sizes $N_0$, $N_A$ and $N_B$, let $\bar c$ be the mean of the first
part and put $Z_i=Y_i-\bar c$ in the other two, so that $Z=\varepsilon-d$
with $d=\bar c-\mu_{\mathrm{ref}}$. Let $K_8>0$ satisfy
$\E\varepsilon^8\leq K_8$. On the first part compute a
two-sided Catoni interval $[L_0,U_0]$ for $\mu_{\mathrm{ref}}$ with
$v=K_8^{1/4}$ and put $\bar d=\max(|\bar c-L_0|,|\bar c-U_0|)$. Set
$K_8'=(K_8^{1/8}+\bar d)^8$ and, on the second part, compute an upper Catoni
bound $\widehat m_4^{A}$ for $\E Z^4$ with $v=K_8'$. On the third part compute
two-sided Catoni intervals $J_1,J_2,J_3$ for the raw moments $m_k=\E Z^k$
with $v=(\widehat m_4^{A})^{1/2}$, $\widehat m_4^{A}$ and
$(\widehat m_4^{A}K_8')^{1/2}$, and an upper bound $U$ for $m_4$ with
$v=K_8'$. Let $[c-r,c+r]$, $[c_3-r_3,c_3+r_3]$ and $(-\infty,U_4]$ enclose,
by interval arithmetic over $J_1\times J_2\times J_3\times[0,U]$, the values of
\[
 \sigma^2=m_2-m_1^2,\qquad
 \mu_3=m_3-3m_1m_2+2m_1^3,\qquad
 \mu_4=m_4-4m_1m_3+6m_1^2m_2-3m_1^4 .
\]
If the six error levels sum to $\delta$, then with probability at least
$1-\delta$ the true $(\sigma^2,\mu_3,\mu_4)$ lies in this region.
\end{lemma}

\begin{proof}
By Lyapunov's inequality
$\Var Y=\sigma^2\leq(\E\varepsilon^8)^{1/4}\leq K_8^{1/4}$, which validates
the interval for $\mu_{\mathrm{ref}}$; on its event $|d|\leq\bar d$. Given the
first part, $\bar c$, $\bar d$ and $K_8'$ are fixed and the $Z_i$ in the other
two parts are independent copies of $\varepsilon-d$. On the first event,
Minkowski's inequality gives $(\E Z^8)^{1/8}\leq K_8^{1/8}+\bar d$, so
$\Var Z^4\leq K_8'$ with a proxy fixed before the second part is used, and
the crude bound holds conditionally at its level. On both events $\Var Z\leq(\E Z^4)^{1/2}$,
$\Var Z^2\leq\E Z^4$, $\Var Z^3\leq\E Z^6\leq(\E Z^4\E Z^8)^{1/2}$ and
$\Var Z^4\leq\E Z^8$, and these proxies are fixed given the first two parts,
so the bounds on the third part hold conditionally at their levels. The
displayed identities expand the central moments of $Z$ about its mean $m_1$,
and these equal the central moments of $\varepsilon$. Interval arithmetic
encloses the image of the box, and a union bound finishes the proof.
\end{proof}

The centering costs little. The interval for $m_1$ is as narrow as a
confidence interval for a mean with variance $\sigma^2$, and its half-width
enters $r_3$ mainly through the term $3m_1m_2$, an increase of order
$3(\sigma^2/N_B)^{1/2}$ against an $r_3$ of order
$\{(\mu_4K_8)^{1/2}/N_B\}^{1/2}$, up to the same logarithmic factor.
The first-order comparator uses $\sigma^2\leq m_2$ from the same chain.
Everything below uses only the region, so it applies to either lemma.

On the event of Lemma~\ref{lem:region}, use
\[
 \phi_{m,c}(X)=\begin{pmatrix}X-m\\(X-m)^2-c\end{pmatrix},
 \qquad
 \E\phi_{m_0,c}=(0,\Delta)^\top,\quad \Delta=\sigma^2-c,\quad|\Delta|\leq r.
\]
Its second-moment matrix is
$\bigl(\begin{smallmatrix}\sigma^2&\mu_3\\ \mu_3&\mu_4-\sigma^4+\Delta^2
\end{smallmatrix}\bigr)$. With $s_+=c+r$, $s_-=\max(c-r,0)$ and
$T=U_4-s_-^2+r^2$, the matrix
\[
 \Fbar=\begin{pmatrix} s_++a & c_3\\ c_3 & T+r_3^2/a\end{pmatrix},
 \qquad a>0,
\]
dominates it for every admissible moment triple, because the difference has
nonnegative diagonal entries whose product is at least $r_3^2$. We take
$a=\max\{r_3(|c_3|-r_3)/T,\ 10^{-9}s_+\}$, which minimizes
$\Fbar_{11}-\Fbar_{12}^2/\Fbar_{22}$; when $|c_3|\leq r_3$ the score
effectively reduces to the first-order one. Proposition~\ref{prop:robust}
with $\mathcal B=\{(0,y):|y|\leq r\}$ and $\rho_t=|\lambda_{t,2}|r$ then
gives an exact e-process, and Corollary~\ref{cor:calibration}'s union bound
gives coverage at least $1-\delta-\alpha_e$ for the inverted confidence
sequence at level $\alpha_e$. No profiling over $\eta$ is needed.

The price of calibration is a loss of local width, and it has a closed form.
Write $\Fbar=(f_{ij})$,
\[
 \bar{\mathrm g}=1-\frac{f_{12}^2}{f_{11}f_{22}},
 \qquad
 u(\rho)=
 \begin{cases}
 (1-\bar{\mathrm g})\rho+\{\bar{\mathrm g}\,(1-(1-\bar{\mathrm g})\rho^2)\}^{1/2},
   & 0\leq\rho\leq1,\\
 1,&\rho>1.
 \end{cases}
\]

\begin{proposition}[Calibrated local width]\label{prop:calwidth}
Let $\Fbar$ be symmetric positive definite with $f_{12}\neq0$, let
$|\Delta|\leq r$, fix a shift $w>0$ and let $I=[\Delta-r,\Delta+r]$. Over
directions $a\in\R^2$, the leading local growth
$\sup_{\kappa>0}\{\kappa\min_{y\in I}a^\top(w,y)-\kappa^2a^\top\Fbar a/2\}$
is maximized at value
$\tfrac12\min_{y\in I}(w,y)\Fbar^{-1}(w,y)^\top$, and
\[
 (w,y)\Fbar^{-1}(w,y)^\top=\frac{w^2}{f_{11}}
 +\frac{(y-f_{12}w/f_{11})^2}{f_{22}-f_{12}^2/f_{11}}.
\]
Consequently, if the boundary sequence is matched with the first-order
comparator as in Proposition~\ref{prop:width}, so that the latter has local
half-width $\{2b_nf_{11}/n\}^{1/2}$, the best calibrated half-width divided by
it equals $u(\rho_n)$, where
\[
 \rho_n=\frac{r_+\,f_{11}^{1/2}}{|f_{12}|\,(2b_n/n)^{1/2}},
 \qquad r_+=r+\operatorname{sign}(f_{12})\Delta\in[0,2r].
\]
In particular $u(0)=\bar{\mathrm g}^{1/2}$ and $u(\rho)=1$ for $\rho\geq1$.
The half-width on the side $w<0$ is given by the same formula with $\Delta$
replaced by $-\Delta$.
\end{proposition}

\begin{proof}
For any $a$ and $y\in I$, the Cauchy--Schwarz inequality gives
$a^\top(w,y)\leq(a^\top\Fbar a)^{1/2}\|(w,y)\|_{\Fbar^{-1}}$, so the growth
is at most $\tfrac12\min_{y\in I}\|(w,y)\|^2_{\Fbar^{-1}}$. Let $y^*$
attain this minimum over the convex set $I$ and set
$a^*=\Fbar^{-1}(w,y^*)$. The first-order optimality condition gives
$a^{*\top}(w,y)\geq a^{*\top}(w,y^*)=\|(w,y^*)\|^2_{\Fbar^{-1}}$ for all
$y\in I$, and $a^{*\top}\Fbar a^*$ equals the same quantity, so $a^*$
attains the bound with $\kappa=1$. The displayed identity is the Schur
complement formula. The minimizing $y$ is the point of $I$ nearest
$f_{12}w/f_{11}$. For $w>0$ this point has the sign of $f_{12}$, so it leaves
$I$ through $\Delta+r$ when $f_{12}>0$ and through $\Delta-r$ when $f_{12}<0$,
and its distance from $I$ is $(|f_{12}|w/f_{11}-r_+)_+$ in both cases. Writing
$w=u\{2b_nf_{11}/n\}^{1/2}$, the boundary equation
$\min_{y\in I}\|(w,y)\|_{\Fbar^{-1}}^2=2b_n/n$ becomes
$\bar{\mathrm g}u^2+(1-\bar{\mathrm g})(u-\rho_n)_+^2=\bar{\mathrm g}$,
whose root in $[\rho_n,1]$ is $u(\rho_n)$; for $\rho_n\geq1$ the root is
$u=1$. The map $(w,y)\mapsto(-w,-y)$ leaves the quadratic form unchanged and
sends $I$ to $[-\Delta-r,-\Delta+r]$, which gives the side $w<0$.
\end{proof}

Proposition~\ref{prop:calwidth} makes the trade-off quantitative. With
$v\approx\mu_4$ in Lemma~\ref{lem:region} the variance half-width is
\[
 r\approx\{2\mu_4\log(2/\delta_\sigma)/N_B\}^{1/2},
\]
so for standardized errors
\[
 \rho_n\approx\frac{1}{|\gamma_3|}
 \Bigl\{\frac{(3+\gamma_4)\log(2/\delta_\sigma)}{b_n}\Bigr\}^{1/2}
 \Bigl(\frac{n}{N_B}\Bigr)^{1/2}.
\]
The reference size therefore enters twice. Through $\bar{\mathrm g}$, the
uncertainty about $\mu_3$ and $\mu_4$ caps the gain at $u(0)=\bar{\mathrm
g}^{1/2}$, which tends to $\gtwo^{1/2}$ only as $N\to\infty$. Through $\rho_n$,
the uncertainty about $\sigma^2$ erodes what remains: since
$u'(0)=1-\bar{\mathrm g}$, the loss grows at the rate $(n/N_B)^{1/2}$, and the
gain disappears once $n$ is of the order of $N_B$. The oracle constant
$\sqrt{\gtwo}$ of Proposition~\ref{prop:width} is the limit $N\to\infty$ at
fixed $n$.

The minimizing direction depends on $w$ and thus on the horizon. In the
implementation, the component with scale $a_H$ of the first-order bank
targets $w_H=a_Hf_{11}^{1/2}u(\rho_H)$, where $\rho_H$ is $\rho_n$ with
$(2b_n/n)^{1/2}$ replaced by $a_H$ and $\Delta=0$. Its stake is
$\lambda_H=\Fbar^{-1}(w_H,y_H)$, with $y_H$ the point of $[-r,r]$ nearest
$f_{12}w_H/f_{11}$. The bank is fixed once the reference sample is
processed, so it is predictable for the tested stream. With $r=0$ and
$\Fbar=F$ it coincides with the oracle bank of
Section~\ref{sec:experiments}. Stakes with $|\lambda_{H,2}|r\geq0.45$ would be
shrunk to keep $\rho_t<1/2$; none required shrinking in the study below,
where the largest penalty $|\lambda_{H,2}|r$ over all reference samples was
0.0027.

\section{Design of the numerical study}\label{sec:experiments}

The numerical study has four parts. Simulations~1--5 check the theory when
the moments are known. Simulations~6 and~7 examine misspecified moments and
alternative betting strategies. Simulations~8--10 replace the known moments by
a reference sample, and a final example applies the calibrated method to
gamma-ray detector data. In every part the laws, sample sizes, random seeds and
acceptance criteria were fixed before the corresponding computations were run.
Simulations~1--5 form the main study. The remaining parts were specified after
it, each before its own runs, and are reported separately.

\subsection{Common settings}

For $H\in\{100,300,1000,3000,10000\}$, define
\[
 a_H=\{2\log(2/\alpha)/H\}^{1/2}.
\]
The two-sided procedure was the equal-weight mixture of the ten strategies
with scales $\pm a_H$. The moment method set
$\lambda=\pm a_HF^{-1}e_1/\sqrt{e_1^\top F^{-1}e_1}$; the first-order
comparator used the same rational envelope with unit-variance score
$(X-m)/\sigma$. Thus every component was fixed before observing its current
factor and met the predictability requirement.

The Catoni comparator used the same scale bank and variance bound, with
log-factor $\psi(a_H(X-m)/\sigma)-a_H^2/2$, where
\[
 \psi(u)=
 \begin{cases}
 \log(1+u+u^2/2),&u\geq0,\\
 -\log(1-u+u^2/2),&u<0.
 \end{cases}
\]
This is the influence function of \citet{wang2023}. The e-processes of
\citet{fan2025} test one-sided mean/variance inequalities and are not direct
competitors for the projected location-CS width studied here; we use them to
delimit prior art rather than force a mismatched numerical comparison.

\subsection{Known moments}
\label{sec:design-oracle}

Simulation~1 assessed validity under the null. For each of five laws
(Gaussian, Gamma(2), lognormal, a skew Gaussian mixture and Pareto(5)) it
generated 3000 null paths through time 10000 and recorded whether the mixture
e-process ever reached $1/\alpha=20$. The acceptance criterion was a one-sided
95\% Clopper--Pearson upper bound of at most 0.07 for the crossing
probability. Pareto(5) rather than Pareto(2.5) was used because the theorem
needs a finite fourth moment.

Simulation~2 measured confidence-sequence widths by inverting terminal mixture
tests on a 161-point grid at $n\in\{100,1000,10000\}$. For each standardized
path, the grid was centered at the sample mean and extended $10/\sqrt n$ in
both directions; the width is the convex hull of nonrejected grid points,
without interpolation. Sixteen paired paths per cell suffice to check the
width constant but not to rank methods in finite samples. The criterion was a
relative error of at most 15\% from $\sqrt{\gtwo}$ on skewed laws and, under
the Gaussian control, a moment construction at most 5\% wider than the
first-order one. The data-centered grid is a diagnostic approximation rather
than the coverage-valid inversion: it could understate the width if a
disconnected nonrejected region lay beyond the grid.

Simulation~3 evaluated the projection identity \eqref{eq:kappa} on three
finite-support exponential paths with cubic, threshold and sinusoidal tangent
directions. Simulation~4 compared the exact expected log-growth with the
approximation \eqref{eq:normal} over
$\delta/\sigma\in\{0.02,0.05,0.1,0.2,0.4,0.7\}$. In this simulation only,
$\lambda=F^{-1}\beta$, so that
$\E_Qh=\lambda^\top\beta=\lambda^\top F\lambda=v$; this identity served as a
control variate, which changes the Monte Carlo variance but not the estimand.

Simulation~5 compared stopping times on two standardized
zero-inflated-lognormal laws, with nonzero probabilities 0.25 and 0.60, under
shifts of 0.08, 0.12 and 0.18 standard deviations. Each cell contained 2000
paired paths through time 10000. The criterion was a median reduction of the
stopping time of at least 15\% in at least one scenario, with all scenarios
reported.

\subsection{Misspecified moments and betting strategies}
\label{sec:design-diagnostics}

Two supplementary simulations, specified after the main study, examine
settings that the theorem does not cover.

Simulation~6 used 600 paths through time 10000 for Gaussian, Gamma(2), and
lognormal null laws. Variance, third moment, and fourth moment were perturbed
one at a time by $\pm2\%$, $\pm5\%$, and $\pm10\%$. For the Gaussian
third moment, which equals zero, the perturbation was scaled by the natural
standardized covariance scale $\sqrt{\mu_4-1}$. Separate runs inflated the
true matrix by 1.05, 1.10, and 1.25 while preserving exact centering. The
outputs include the ever-crossing rate, a one-sided Clopper--Pearson bound,
the two terms in \eqref{eq:misspec}, and stopping times under a shift of
$0.12\sigma$.

Simulation~7 compared the 10-component bank of the main study, a
30-component bank, and quadratic follow-the-regularized-leader (FTRL) updates
with regularization 10, 50 or 200 and stake cap 0.35. Null runs used the five
laws of Simulation~1, and stopping-time runs the two laws and three shifts of
Simulation~5. Each cell contained 600 paths through time 10000.
Writing $a_{\rm proj}=F^{-1}e_1/\sqrt{e_1^\top F^{-1}e_1}$, the scalar FTRL stake was
\[
 \lambda_t=\operatorname{clip}\left\{
 \frac{\sum_{s<t}z_s}{r+t-1},-0.35,0.35\right\},
\]
where
$z_s=e_1^\top F^{-1}\phi_s/\sqrt{e_1^\top F^{-1}e_1}$ is the
unit-second-moment projected score and $r$ is the regularization. This formula
sets the vector stake to $\lambda_t a_{\rm proj}$, so the current score
does not enter the current stake.

\subsection{Moments from a reference sample}
\label{sec:design-calibrated}

Simulations~8--10 replace the known moments by the construction of
Section~\ref{sec:calibrated}. Each replication drew its own reference sample of
size $N\in\{10^4,10^5,10^6\}$ from the standardized error law, with
$N_A=0.2N$ and an envelope $K_8$ equal to twice the true eighth central
moment. The total calibration error $\delta=0.01$ was split as $0.002$ for
$\widehat\mu_4^{A}$, $0.004$ for $\sigma^2$, and $0.002$ each for $\mu_3$ and
$\mu_4$. The e-processes used $\alpha_e=0.04$, so the overall level was 0.05.
The comparator was the first-order factor with unit-variance score
$(X-m)/\bar s^{1/2}$, where the upper variance bound $\bar s$ came from the
same reference sample through the same first step ($0.002$) and a one-sided
bound at level $0.008$, so that both methods spent the same calibration
budget. Both used the scale bank $a_H$ with $\alpha_e$ in place of $\alpha$.
Versions with the true moments, at level $\alpha_e$, were run on the same
streams as a benchmark.

The laws were Gamma(2), the lognormal law, and the two zero-inflated
lognormal laws of Simulation~5. Their standardized eighth central moments are
$3.75\times10^3$, $8.14\times10^3$, $3.19\times10^6$ (nonzero probability
0.60), and $2.49\times10^{12}$ (nonzero probability 0.25). Simulation~8 used
1000 replications per law and $N$. For validity it recorded whether either
e-process ever reached $1/\alpha_e$ through time 10000 at the true location,
with the criterion of Simulation~1. For power it recorded one-sided stopping
times on fresh streams from the same reference samples, shifted by $0.08$,
$0.12$ and $0.18$. The criterion was a median gain of at least 15\% over the
calibrated first-order comparator at $N=10^6$ for Gamma(2) or the lognormal
law.

Simulation~9 repeated the width measurement of Simulation~2 with 32
replications per law and $N$, excluding the law with nonzero probability
0.25. It compared the mean width ratio with a matched-boundary prediction from
Lemma~\ref{lem:remainder} without its remainder. The boundary level was read
off the observed first-order half-width, and the calibrated half-width was the
smallest shift at which the best bank component reached that level. The
criterion was 15\% relative agreement for $n\geq1000$.

Simulation~10 repeated Simulation~8 with the reference record of
Lemma~\ref{lem:region-unknown}, on the same random draws, so that each
replication saw the same reference errors and the same streams. The reference
was observed with an offset of 3.7 that the procedure did not know. The splits
were $N_0=0.02N$, $N_A=0.2N$ and the rest, and the calibration error
$\delta=0.01$ was split as $0.001$ each for the centering interval and
$\widehat m_4^{A}$, and $0.002$ each for $m_1$, $m_2$, $m_3$ and $m_4$. The
first-order comparator used the same first two steps and an upper bound for
$m_2$ at level $0.008$. The criteria were those of Simulation~8. A first run
had bounded $|d|$ on the second part, which left the crude bound
$\widehat m_4^{A}$ without a fixed variance proxy. The bound was moved to the
first part as in Lemma~\ref{lem:region-unknown}, and Simulation~10 and the
real-data example were repeated on the same draws and files with no other
change.

\subsection{Gamma-ray detector data}
\label{sec:design-real}

The real-data example uses public data of the Fermi Gamma-ray Burst Monitor
\citep{meegan2009gbm}, distributed by the HEASARC archive of NASA. Each photon
registered by the sodium iodide detector n0 carries a pulse-height channel
between 0 and 127 that increases with its energy; outside bursts the channel
distribution is right-skewed. Consecutive channels were the observations. The
design was fixed on a development day, 15 March 2024, and then applied
unchanged to a confirmation day, 15 June 2024, whose data were not examined
beforehand. The first hour of the confirmation day fixed the unit of
measurement (its mean and standard deviation) and the envelope $K_8$, twice
its standardized eighth central moment, and supplied no windows; its skewness
and kurtosis give $\gtwo=0.39$. The remaining hours were cut into consecutive
windows, each consisting of a 300-second reference record, about
$3\times10^5$ events, followed by a monitored stream of $n=10^4$ events.
Windows with a low event rate (passages through the South Atlantic Anomaly)
or a data gap were excluded. Calibration followed
Lemma~\ref{lem:region-unknown} with the settings of Simulation~10. Because
the stream may differ slightly from its reference, a stability margin,
declared in advance, enlarged the region to contain the central moments of
each of ten consecutive segments of the reference, and the first-order
variance bound to the largest segment variance.

On the development day the mean of a stream differed from the mean of the
preceding 300 seconds by far more than noise allows (standard deviation 0.031
against 0.010). The location of the channel distribution drifts within
minutes, so the reference mean is not the location of the stream. Validity
was therefore assessed at the mean of all events within 30 seconds before and
after the stream, which follows the current location, with the criterion of
Simulation~1 applied to the calibrated moment e-process with margin. The
comparator, the variant without margin, and the reference mean as null value
are reported for comparison. Power was assessed by adding shifts of 0.08,
0.12 and 0.18 to the same streams, with the criterion of a median gain of at
least 15\% at 0.12.

\section{Results}\label{sec:results}

All acceptance criteria of Section~\ref{sec:experiments} were met. The
subsections follow the four parts of the design.

\subsection{Known moments}
\label{sec:results-oracle}

\begin{figure}[t]
 \centering
 \includegraphics[width=0.78\textwidth]{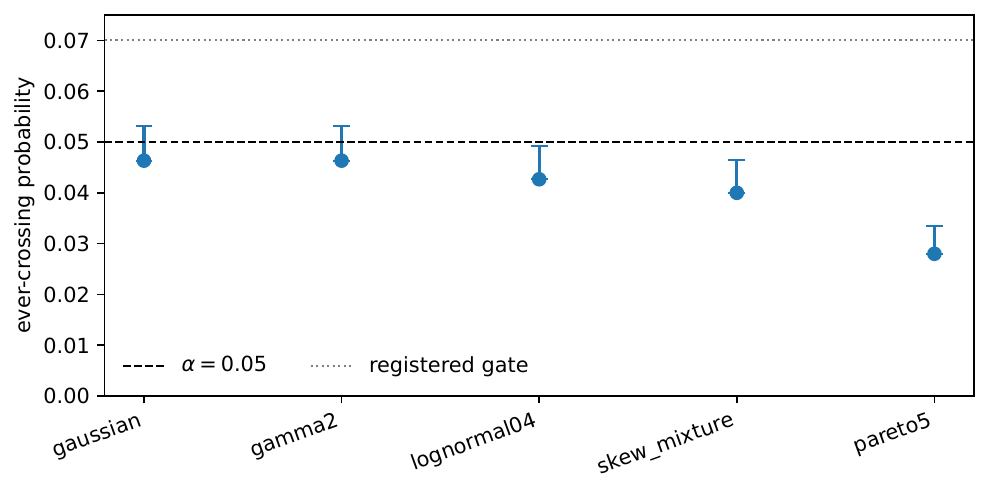}
 \caption{Estimated ever-crossing probabilities for Gaussian, Gamma(2),
 lognormal, skew Gaussian-mixture, and Pareto(5) null laws. Each point uses
 3000 paths through time 10000 at $\alpha=0.05$; bars end at the one-sided
 95\% Clopper--Pearson upper bound. The horizontal lines mark the nominal
 level and the acceptance bound.}
 \label{fig:validity}
\end{figure}

The simulated crossing rates under the null were all below the nominal 5\%
(Simulation~1, Figure~\ref{fig:validity}). They ranged from 0.0280 to 0.0463,
and the largest one-sided 95\% Clopper--Pearson upper bound was 0.0531, below
the acceptance bound of 0.07. A finite simulation of this kind cannot confirm
the guarantee, which follows from Theorem~\ref{thm:eprocess}; it checks the
implementation.

\begin{figure}[t]
 \centering
 \includegraphics[width=0.78\textwidth]{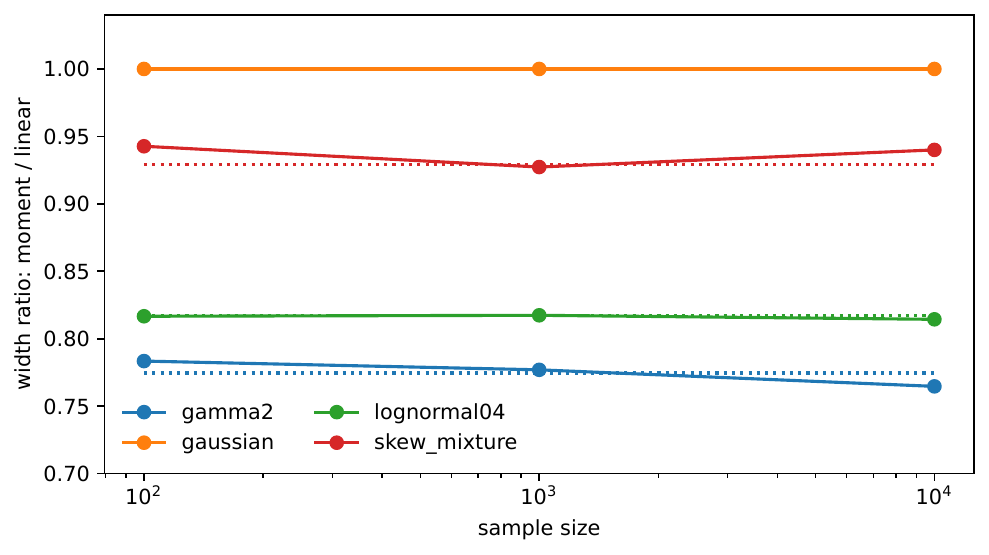}
 \caption{Empirical ratios of moment to first-order confidence-sequence width
 track the theoretical $\sqrt{\gtwo}$ constants, the predicted asymptotic gain
 from moment orthogonalization. Solid lines are paired Monte Carlo means from
 16 paths per law-by-sample-size cell; dotted lines are the limits in
 Proposition~\ref{prop:width}.}
 \label{fig:width}
\end{figure}

The width ratios followed Proposition~\ref{prop:width} closely
(Simulation~2, Figure~\ref{fig:width}). At $n=10000$ the observed
moment-to-linear ratios were 1.000 for the Gaussian control, 0.765 for
Gamma(2), 0.814 for the lognormal law and 0.940 for the skew Gaussian
mixture, against limits of 1.000, 0.775, 0.817 and 0.929. No accepted set
touched the boundary of the location grid.

The projection identity held to floating-point precision (Simulation~3). The
captured fractions were 0.722676 for the cubic five-point path, 0.649927 for
the threshold five-point path and 0.592589 for the sinusoidal six-point path,
and the direct projection deficit and $(1-\kappa)\chi^2/2$ differed by at most
$4.45\times10^{-17}$.

The quadratic growth approximation was accurate over a nontrivial range of
shifts (Simulation~4). Its relative error was at most 2.0\% for every law and
$|\delta|\leq0.1\sigma$. It first exceeded 10\% at $\delta=0.4$ for the
Gaussian law and at $\delta=0.7$ for the skew mixture, and stayed below 10\%
on the whole grid for Gamma(2) and the lognormal law. The analytic remainder
bound held in every cell.

With known moments, the moment score stopped earlier in all six
law-by-shift cells of Simulation~5 (Table~\ref{tab:ab}). Median reductions
relative to the first-order rational factor ranged from 36.8\% to 45.8\%,
and paired bootstrap intervals excluded zero in every cell. All methods used
identical paired streams and the same scale bank. The Catoni and first-order
rational factors, which share the variance bound, differed by at most 6.5
observations in median stopping time.

\begin{table}[t]
\centering
\caption{Median stopping times for 2000 paired zero-inflated-lognormal paths.
Smaller is better. Gain is $1-T_{\rm moment}/T_{\rm linear}$; brackets give a
paired-bootstrap 95\% interval based on 2000 resamples.}
\label{tab:ab}
\small
\begin{tabular}{lrrrrr}
\toprule
Nonzero probability & $\delta$ & Moment & Linear & Catoni & Gain [95\% CI] \\
\midrule
0.25 & 0.08 & 750.5 & 1203.0 & 1205.0 & 37.6 [35.7,39.7]\% \\
0.25 & 0.12 & 362.0 & 573.0 & 575.0 & 36.8 [35.0,39.1]\% \\
0.25 & 0.18 & 164.0 & 265.0 & 267.0 & 38.1 [35.9,39.9]\% \\
0.60 & 0.08 & 565.5 & 1043.0 & 1049.5 & 45.8 [42.9,47.9]\% \\
0.60 & 0.12 & 270.5 & 491.0 & 494.5 & 44.9 [42.6,47.0]\% \\
0.60 & 0.18 & 133.0 & 218.0 & 221.5 & 39.0 [37.5,41.5]\% \\
\bottomrule
\end{tabular}
\end{table}

\subsection{Misspecified moments and betting strategies}
\label{sec:results-diagnostics}

Equation~\eqref{eq:misspec} predicted the qualitative behavior under
misspecified moments (Simulation~6). Runs with the true moments had zero
centering and matrix errors, one-step
conditional mean exactly one, and ever-crossing estimates between 0.0350 and
0.0467. Perturbing the variance changed both centering and the matrix. With a
$-2\%$ error, the Gamma(2) crossing estimate rose to 0.0800 (one-sided 95\% upper
bound 0.101); a $-10\%$ variance error under Gamma(2) produced the largest
estimate, 0.840. These are
finite-horizon diagnostics of invalidity under misspecification, not revised
validity guarantees.

Third- and fourth-moment perturbations left the score centered but could
violate the PSD bound. Across the three laws, the largest crossing estimates
at $\pm10\%$ were 0.083 for the third moment and 0.072 for the fourth.
Conversely, multiplying the true matrix by 1.05, 1.10, and 1.25 retained the
PSD inequality. Their crossing ranges were 0.027--0.038, 0.025--0.030, and
0.007--0.015, respectively. The corresponding median stopping-time penalties
relative to the oracle ranged over 5.0--7.4\%, 9.3--13.2\%, and
31.1--33.2\%. This is the expected validity--power tradeoff from conservative
inflation.

All null estimates in the strategy comparison (Simulation~7) lay between
0.0200 and 0.0533. The largest
one-sided 95\% upper bound was 0.0710, for FTRL with regularization 10 under
the lognormal null. With only 600 paths per cell, this diagnostic retains
nontrivial Monte Carlo uncertainty and is not a replacement for the theorem.
Table~\ref{tab:strategy} shows that the denser bank
changed median stopping by $-2.1\%$ to $5.5\%$ relative to the 10-component
bank. FTRL was more sensitive: regularization 50 improved the median by
1.5--15.6\% across the six cells, whereas regularization 200 was slower in
four cells and lost as much as 38.0\%. The theorem establishes validity for
predictable adaptation; these finite-horizon diagnostics only audit its
implementation and show that it does not remove the need to tune or mix
strategies.

\begin{table}[t]
\centering
\caption{Median stopping times in the strategy comparison of Simulation~7
(600 paired paths per cell). ZILN($p$) denotes the standardized
zero-inflated-lognormal law with nonzero probability $p$. Follow-the-regularized-
leader (FTRL) uses the displayed prefix-only update; its three columns use
regularization 10, 50, and 200. Fixed 10 is the scale bank of the main
study, and Dense 30 is the corresponding 30-component bank.}
\label{tab:strategy}
\small
\begin{tabular}{lrrrrrr}
\toprule
Law & $\delta$ & Fixed 10 & Dense 30 & FTRL-10 & FTRL-50 & FTRL-200 \\
\midrule
ZILN(0.25) & 0.08 & 930.5 & 892.0 & 912.5 & 813.0 & 844.5 \\
ZILN(0.25) & 0.12 & 429.5 & 406.0 & 375.5 & 362.5 & 444.5 \\
ZILN(0.25) & 0.18 & 203.0 & 203.0 & 174.5 & 188.0 & 250.0 \\
ZILN(0.60) & 0.08 & 762.5 & 743.0 & 752.5 & 703.5 & 763.0 \\
ZILN(0.60) & 0.12 & 326.0 & 317.5 & 309.0 & 308.0 & 382.5 \\
ZILN(0.60) & 0.18 & 163.0 & 166.5 & 144.0 & 160.5 & 225.0 \\
\bottomrule
\end{tabular}
\end{table}

A small timing check used four Gamma(2) paths, 81 candidate means,
and 1000 observations. The measured costs were 0.077, 0.209, and
0.034--0.040 milliseconds per candidate for the 10-component bank,
30-component bank, and each FTRL update, respectively. These timings are implementation-
and machine-dependent; this limited benchmark showed no obvious computational
bottleneck from predictability.

\subsection{Moments from a reference sample}
\label{sec:results-calibrated}

Calibration preserved validity (Simulation~8). The region of
Lemma~\ref{lem:region} contained the true moment triple in all 12000
replications, and the first-order variance bound held in all of them. The
calibrated moment e-process reached $1/\alpha_e$ at the true location in at
most 1.8\% of replications per cell (largest one-sided 95\%
Clopper--Pearson bound 0.027), and the calibrated first-order comparator in
at most 2.9\% (bound 0.039). Both rates are below $\alpha_e=0.04$, as
expected when the bounds dominate the true moments.

\begin{figure}[t]
 \centering
 \includegraphics[width=\textwidth]{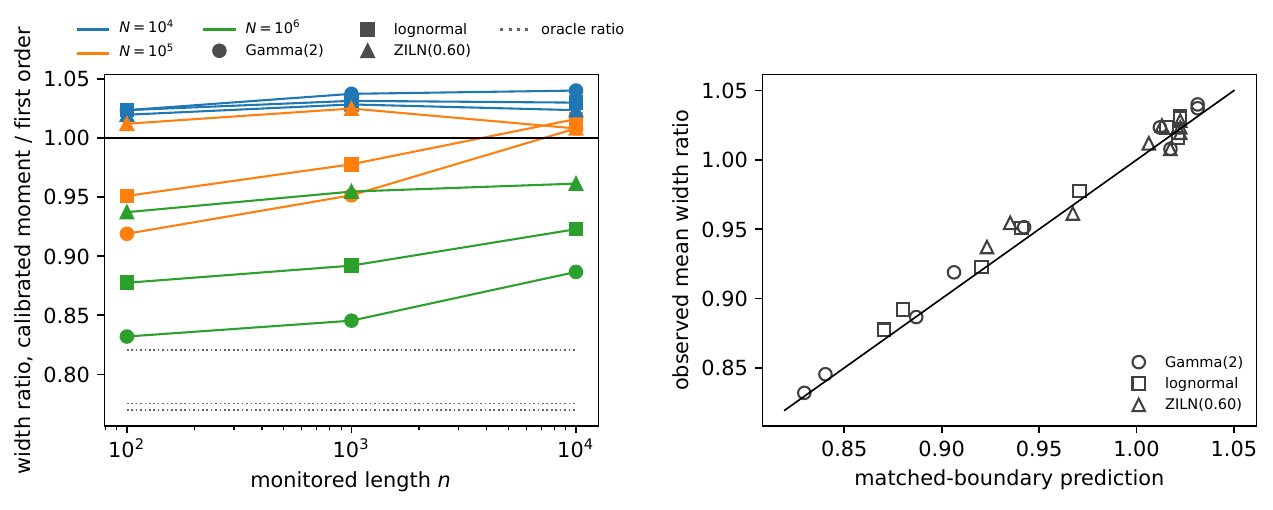}
 \caption{Simulation~9. Left: mean ratio of calibrated moment to calibrated
 first-order confidence-sequence width against the monitored length, for three
 reference sizes $N$ (colors) and three error laws (markers); dotted lines are
 the oracle ratios on the same streams. Right: observed mean ratio against the
 matched-boundary prediction for all 27 cells; the line is the identity. Each
 cell uses 32 replications, each with its own reference sample.}
 \label{fig:cal}
\end{figure}

The width law of Proposition~\ref{prop:calwidth} was confirmed
(Simulation~9, Figure~\ref{fig:cal}). In all 27 cells, the
mean ratio of the calibrated moment width to the calibrated first-order width
was within 2.1\% of the matched-boundary prediction, against a 15\%
tolerance, and no accepted set touched the grid edge. The closed form of
Proposition~\ref{prop:calwidth}, evaluated at the design radius, differed
from the component-level prediction by at most 0.007. The observed ratio
ranged from 0.832 (Gamma(2), $N=10^6$, $n=100$) to 1.040 (Gamma(2),
$N=10^4$, $n=10^4$); the oracle ratio on the same streams was 0.747--0.830.
The reference size acts in two ways. Uncertainty about $\mu_3$ and $\mu_4$
raises $\bar{\mathrm g}$ toward one, which caps the gain at every $n$; over
the reference samples of Simulation~8 with $N=10^4$, the median $\bar{\mathrm g}$ was
0.91 for Gamma(2) against
$\gtwo=0.60$. Uncertainty about $\sigma^2$ then erodes the remaining gain as
$n$ grows, in line with the $(n/N)^{1/2}$ dependence of $\rho_n$.

\begin{table}[t]
\centering
\caption{Simulation~8: median stopping times at shift $0.12$ (1000 replications per
cell, each with its own reference sample of size $N$). ``Calibrated'' uses
the centering-robust factor with moments from the reference sample; ``first
order'' uses the variance bound from the same sample and budget. Oracle
columns use the true moments on the same streams (oracle). Gains are
$1-T_{\rm moment}/T_{\rm first\ order}$. A median of 10001 means that most
paths did not stop by time 10000.}
\label{tab:cal}
\small
\begin{tabular}{llrrrrrr}
\toprule
 & & \multicolumn{3}{c}{calibrated} & \multicolumn{3}{c}{oracle} \\
\cmidrule(lr){3-5}\cmidrule(lr){6-8}
Law & $N$ & moment & first order & gain & moment & first order & gain\\
\midrule
Gamma(2) & $10^{4}$ & 616.5 & 567.0 & $-8.7\%$ & 306.0 & 505.5 & 39.5\% \\
Gamma(2) & $10^{5}$ & 448.0 & 518.5 & 13.6\% & 288.0 & 494.5 & 41.8\% \\
Gamma(2) & $10^{6}$ & 355.0 & 483.0 & 26.5\% & 292.5 & 477.0 & 38.7\% \\
lognormal & $10^{4}$ & 594.5 & 545.0 & $-9.1\%$ & 303.0 & 472.0 & 35.8\% \\
lognormal & $10^{5}$ & 478.5 & 505.0 & 5.2\% & 315.5 & 472.5 & 33.2\% \\
lognormal & $10^{6}$ & 388.5 & 498.5 & 22.1\% & 328.0 & 488.0 & 32.8\% \\
ZILN(0.60) & $10^{4}$ & 853.5 & 826.0 & $-3.3\%$ & 291.0 & 519.0 & 43.9\% \\
ZILN(0.60) & $10^{5}$ & 599.0 & 593.5 & $-0.9\%$ & 289.5 & 535.0 & 45.9\% \\
ZILN(0.60) & $10^{6}$ & 461.0 & 534.0 & 13.7\% & 298.0 & 518.0 & 42.5\% \\
ZILN(0.25) & $10^{4}$ & 10001.0 & 10001.0 & 0.0\% & 389.0 & 621.0 & 37.4\% \\
ZILN(0.25) & $10^{5}$ & 2331.5 & 2089.5 & $-11.6\%$ & 368.0 & 587.5 & 37.4\% \\
ZILN(0.25) & $10^{6}$ & 896.5 & 859.0 & $-4.4\%$ & 359.0 & 588.0 & 38.9\% \\
\bottomrule
\end{tabular}
\end{table}

Calibration costs part of the gain, and how much depends on $N$
(Table~\ref{tab:cal}). With the true moments the gains were 30.2--47.5\% in
every cell. Across the three shifts they shrank
to 26.5--30.3\% for
Gamma(2), 21.6--22.2\% for the lognormal law, and 8.5--15.1\% for ZILN(0.60)
at $N=10^6$. At $N=10^5$ they were 5.4--13.6\%, 3.9--9.0\%, and between
$-5.1\%$ and $-0.9\%$. At $N=10^4$ the calibrated moment method was 3.3--9.9\%
slower than the calibrated first-order comparator for these three laws: the
wider two-sided variance interval and the matrix increments then cost more
than the residual skewness information returns. ZILN(0.25), whose eighth
moment is $2.5\times10^{12}$, could not be calibrated: the region was
uninformative ($\bar{\mathrm g}=1.000$ to three decimals at every $N$). The
moment method then degenerated to a slightly less efficient first-order method
and was 3.5--13.4\% slower. At $N=10^4$ neither calibrated method stopped by
time 10000 in most paths at shifts 0.08 and 0.12, because the variance
interval itself was uninformative (median half-width 16.5 for a unit
variance). The
gains of Simulation~5 on this law therefore depend on known moments that a
finite reference sample under an eighth-moment envelope cannot supply.

A reference of unknown location cost little (Simulation~10). The region of
Lemma~\ref{lem:region-unknown} contained the true moment triple in all 12000
replications. The calibrated
moment e-process reached $1/\alpha_e$ at the true location in at most 1.6\% of
replications per cell (largest Clopper--Pearson bound 0.024), and the
first-order comparator in at most 2.9\% (bound 0.039). Not knowing the
reference location cost 1.0--3.0 percentage points of gain at $N=10^6$. Across
the three shifts the gains were 25.5--27.7\% for Gamma(2), 19.0--20.3\% for the
lognormal law and 6.6--12.0\% for ZILN(0.60), against 26.5--30.3\%,
21.6--22.2\% and 8.5--15.1\% with a known location. At $N=10^5$ the loss was
5.1--6.9 points for Gamma(2) and the lognormal law (gains $-0.2$ to 6.9\% and
$-2.4$ to 2.1\%). At $N=10^4$ the median $\bar{\mathrm g}$ was 1.000 to three
decimals for every law, and the moment method was 2.5--11.8\% slower than the
comparator for the first three laws. For Gamma(2) and the lognormal law this
was no worse than with a known location (within 0.1 points), because the
uninformative region makes the method close to the first-order one.
ZILN(0.25) remained uninformative at every $N$.

\subsection{Gamma-ray detector data}
\label{sec:real-results}

On the confirmation day 225 of 234 windows were retained
(Table~\ref{tab:real}). With the stability margin, the calibrated
moment e-process crossed $1/\alpha_e$ at the bracket mean in 4 windows
(1.8\%, Clopper--Pearson bound 0.040) and the first-order comparator in 10
(4.4\%, bound 0.074). Without the margin the moment method crossed in 12
windows (5.3\%, bound 0.085), so the margin was needed. The bracket mean is
itself only a proxy: the standard deviation of the stream mean minus the
bracket mean was 0.0123, against 0.0108 from noise alone, and the residual
drift contributes to both crossing counts. At the reference mean the two
methods with margin crossed in 42 and 52 windows: on these data the reference supplies the
shape of the error law, not the location of the stream. With the margin the
median stopping time fell by 30.3--34.9\% relative to the comparator across
the three shifts, against 31.7--38.0\% without it. The median $\bar{\mathrm g}$
was 0.50, against 0.39 from the characterization hour, and the median
variance half-width was 2.1\% of the variance.

\begin{table}[t]
\centering
\caption{Gamma-ray detector data, confirmation day (225 windows of a
300-second reference record and a stream of $10^4$ photons). Crossings count windows in which the
two-sided e-process reached $1/\alpha_e=25$ within the stream; the bracket
mean is the mean of the events within 30 seconds of the stream. Stopping times
are one-sided, at shift $0.12$ in units of the characterization standard
deviation.}
\label{tab:real}
\small
\begin{tabular}{lrrrrrrr}
\toprule
 & \multicolumn{2}{c}{bracket mean} & \multicolumn{2}{c}{reference mean} &
 \multicolumn{3}{c}{median stop, shift 0.12} \\
\cmidrule(lr){2-3}\cmidrule(lr){4-5}\cmidrule(lr){6-8}
Region & moment & first order & moment & first order & moment & first order & gain \\
\midrule
calibrated & 12 & 11 & 35 & 52 & 353 & 569 & 38.0\% \\
with margin & 4 & 10 & 42 & 52 & 387 & 576 & 32.8\% \\
\bottomrule
\end{tabular}
\end{table}

\section{Discussion}\label{sec:discussion}

The central result is finite-sample. A centered vector score and a
second-moment matrix bound turn the quadratic envelope into an e-factor. This
removes the no-bankruptcy obstruction that affects a raw polynomial bet on an
unbounded support. Moment orthogonalization leaves first-order location
sensitivity unchanged while reducing its variance by $\gtwo$; the confidence
width consequently falls by $\sqrt{\gtwo}$.

Proposition~\ref{prop:quadratic-family} clarifies the relevant optimality
claim. The coefficient $c=1/2$ is locally equivalent to every admissible
quadratic coefficient under the stated local path. Its advantages here are
strict positivity and a short remainder bound, not global optimality. The
strategy comparison makes a parallel distinction: predictability preserves
validity, while its finite-sample power depends on regularization.

The comparison must be interpreted at the level of null information. The
Catoni and first-order methods need a variance bound. The two-component score
also needs valid third- and fourth-moment information, because those moments
center and normalize the additional score. Simulation~5 therefore estimates
the value of this side information. It does not show that higher-order betting can
improve a distribution-free confidence sequence for the mean. Indeed, the
two-point-law characterization in the Introduction rules out that claim.

Simulations~8--10 price that side information when it has to be learned. With an
independent reference sample and an eighth-moment envelope, the exact claim
survives, but its gain is governed by the reference size: $N$ must be large
enough to pin down $\mu_3$ and $\mu_4$, and much larger than the monitored
length $n$. This suits settings in which
reference data are plentiful and decisions are short, such as an instrument
characterized on a standard and then used for brief sequential measurements.
Lemma~\ref{lem:region-unknown} widens this to any long record of the process
in a known-good state, as in the first phase of statistical process control,
at a cost of at most 3 points of gain at $N=10^6$ in Simulation~10. It does not suit errors
whose eighth moment is huge, where even $10^6$ reference observations leave
the moment region uninformative.

The gamma-ray example shows what this means on real data. A skewed law with a
moderate eighth moment and a reference of $3\times10^5$ events gave a gain
close to the simulated one. It also exposed two places where a real reference
departs from the model. The location drifted within minutes, so the reference
supplied the shape of the law but not the null value. And without a stability
margin the moment method crossed more often than $\alpha_e$.

\section{Limitations}

Moment misspecification is the main practical risk. If the second component is
not centered or $F$ fails to dominate its conditional second moment, the
supermartingale proof no longer applies. A predictable upper matrix preserves
validity at the price of power, but estimating the matrix from the same stream
without a valid envelope gives only an asymptotic or calibration-dependent
procedure. The predictive plug-in variant examined during development is
therefore not part of the exact main claim.

Valid moment information can come from an independent historical sample, a
sample-split calibration stage with conservative uncertainty inflation, or
domain constraints that directly imply a predictable $\Fbar_t$. Splitting the
monitored stream itself preserves validity but leaves little gain: with $N$
proportional to $n$, $\rho_n$ in Proposition~\ref{prop:calwidth} is of order
one, and Simulation~9 showed the gain vanishing as $n$ approaches $N$. The gain
requires an external reference with $N\gg n$. The estimation counterpart is
sharper still: when the variance that centers the degree-two score is
estimated from the same residuals, the one-step polynomial estimate of a
location coincides with the sample mean, and sample splitting does not restore
the oracle gain \citep[Corollary~12 and Remark~13]{zabolotnii2026rankone}.
Corollary~\ref{cor:calibration} shows how a valid independent confidence
region transfers to the confidence sequence, and
Lemmas~\ref{lem:region} and~\ref{lem:region-unknown} construct one, the
second without a known reference location. They do not remove the underlying restriction: with unbounded
support and only a finite fourth moment, a finite calibration sample cannot
furnish a nontrivial distribution-free upper bound on that moment, so both
lemmas need the envelope $K_8$. Simulations~8--10 used an envelope twice the true value;
a looser envelope widens the region and moves the break-even point toward
larger $N$. The reference sample must also share the error law of the tested
stream and be independent of it. This requirement tightens as $N$ grows. At
$N=10^6$ the median variance half-width $r$ in Simulation~8 was 1.0--2.4\% of
the variance for the three calibratable laws, so a change of the error law
between reference and stream larger than that is not covered. Enlarging $r$ by a declared stability margin restores validity, and
Proposition~\ref{prop:calwidth} gives the cost.

The simulations have narrower roles than the theorem. Simulation~1 checks the
implementation, Simulation~2 the local width constant, and Simulation~5 gives
a controlled synthetic case in which skewness information matters. Sixteen
paths per cell of Simulation~2 do not support fine finite-sample rankings. The
zero-inflated-lognormal models are standardized synthetic laws rather than
production data. Simulations~6 and~7 were specified after the main study and
have 600, rather than 3000, paths per cell. Simulations~8--10 replace the
known moments by simulated reference samples, so their reference data and
streams are synthetic and independent. The gamma-ray example uses one
detector on one day; its null value is a proxy
for the current location, and its e-processes treat consecutive channels as
independent. Reference data with serial dependence would require a
dependence-robust calibration step that we do not provide.

The relation to prior work is also specific. Catoni's scalar envelope supplies
the positivity mechanism; constrained e-variable theory explains why extra
null information changes the attainable e-power; joint moment e-processes
show that sequential tests under moment constraints are already available;
and Legendre jumpers show that higher-order polynomial betting is established
in exchangeability testing. To our knowledge, the combination that remains
unreported is an exact moment-orthogonal location e-process with the explicit
$\sqrt{\gtwo}$ width calculation. That is the paper's claim boundary.

\section{Extensions}\label{sec:extensions}

The construction may extend beyond location only when the additional scores
can be centered and bounded under the proposed null. Scale inference could use
centered powers or other identification scores if the required higher-moment
bounds are available. For a quantile $q$, the bounded identification score
$\mathbf 1\{X\leq q\}-\tau$ could be combined with jointly centered
auxiliaries. Expectiles and some risk measures offer further candidates
through identification functions. We do not claim explicit width constants,
semiparametric optimality, or automatic orthogonalization for these cases.

\section{Conclusion}

The quadratic envelope converts any correctly centered moment score with a
valid second-moment matrix bound into an exact e-process on an unbounded
support. For the two-component location score, moment orthogonalization changes
the local CS width by the factor $\sqrt{\gtwo}$. The simulated skewed
laws showed the predicted gain, whereas the Gaussian control did not.
Because the variance and the third and fourth central moments belong to the
null specification, the method provides exact sequential inference under that
specification, not a distribution-free dominance over mean-only methods. When the moments come from an independent
reference sample under an eighth-moment envelope, a centering-robust factor
keeps the guarantee exact, and the remaining gain is set by the reference
size and the monitored length through a closed-form width law. The reference
location need not be known, and on real gamma-ray detector data a declared
stability margin kept the error rate below $\alpha_e$ while preserving most of
the gain.

\section*{Code and data availability}

The code is openly available under the MIT licence at\newline
\url{https://github.com/SZabolotnii/Ku-Betting-CS-code-supplement}. It contains the Python implementation, the configurations that
fixed each design before its runs, and the stored results. Scripts there
regenerate every number, table and figure of Sections~\ref{sec:experiments}
and~\ref{sec:results}, and a Lean file checks the algebraic projection and
positivity lemmas. The computations used Python 3.14.5, NumPy 2.4.4 and SciPy
1.17.1, and the Lean file compiles with Lean/Mathlib 4.26.0 without admitted
propositions. The Fermi Gamma-ray Burst Monitor data are public and are
distributed by HEASARC; the code downloads them and records the address and
SHA-256 digest of every file it uses.

\section*{Ethical Considerations}

All numerical examples except the gamma-ray example use synthetic data; that
example uses public instrument data of the Fermi Gamma-ray Burst Monitor. No human-subject,
personal, or production data were used. The method should not be used for high-stakes
sequential decisions unless the stated conditional moment constraints are
independently justified; misspecified centering can invalidate the e-process.
The reported computations were run locally at the scale documented above.

\section*{Declaration of generative AI use}

AI-based tools (Claude Code, Codex) assisted in implementing and checking the
code, in reconciling the numbers of the manuscript with the stored results, in
drafting the text and translating it from Ukrainian, and in language editing
and LaTeX formatting. The author reviewed and edited the whole text and takes
full responsibility for the content of the article.

\section*{Funding}

No external funding was received for this research.

\section*{Disclosure of interest}

The author reports there are no competing interests to declare.

\bibliographystyle{chicago}
\bibliography{references}

\appendix
\section{Measurability of test inversion}\label{app:measurability}

Pointwise coverage requires less regularity than reporting interval endpoints.
Let $\Theta$ be a Borel subset of a Polish space. Suppose
$(\omega,\theta)\mapsto\phi_t(\omega,\theta)$,
$\Fbar_t(\omega,\theta)$, and
$\lambda_t(\omega,\theta)$ are jointly measurable and, for each fixed
$\theta$, $\lambda_t(\theta)$ is predictable. The factor in
\eqref{eq:factor} and every finite
product $M_s(\theta)$ are then jointly measurable. Since time is countable,
$\sup_{s\leq n}M_s(\theta)$ and
$\sup_{s\geq1}M_s(\theta)$ are jointly measurable. Consequently, the
confidence set
\[
 C_n(\omega)=\{\theta:
   \sup_{s\leq n}M_s(\omega,\theta)<1/\alpha\}
\]
has a graph measurable in $\mathcal F\otimes\mathcal B(\Theta)$. These
assumptions suffice for the inversion statement following
Theorem~\ref{thm:eprocess}.

Computable scalar endpoints require additional regularity. Suppose $\Theta$
is a nondegenerate interval in $\R$ with a countable dense subset $D$. Assume
there is a single null set outside which
$\theta\mapsto M_s(\omega,\theta)$ is continuous for every finite $s$. Then
the strict inequality defining $C_n$ and the finite maximum over $s\leq n$
make the nonrejection set relatively open in $\Theta$ at each finite time. Its
infimum and supremum, with fixed conventions for empty or unbounded sets, are
measurable extended-real random variables. For example,
$\{\inf C_n<a\}=\bigcup_{d\in D:d<a}\{d\in C_n\}$ and
$\{\sup C_n>a\}=\bigcup_{d\in D:d>a}\{d\in C_n\}$, with the fixed endpoint
conventions handling the empty set. On the exceptional null set, define both
endpoints by fixed constants; this modification preserves measurability. If an
implementation reports the convex hull of grid points instead, measurability
follows directly because the grid is finite. None of these endpoint
conditions strengthens the null moment assumptions needed for pointwise
e-validity; they are reporting conditions only.

\end{document}